\documentclass[11pt]{article}

\date{September 1st 2020}

\title{\vskip-1.0em\sc Stability of characters and filters for weighted semilattices}
\author{\sc Yemon Choi, Mahya Ghandehari, Hung Le Pham}

\usepackage[a4paper,left=30mm,right=30mm,top=30mm,bottom=35mm]{geometry}

\RequirePackage{amsmath,amssymb, graphicx, mathtools}

\RequirePackage{amsthm}   
\newcounter{pulse}[section]
\numberwithin{pulse}{section}  

\newcommand{\thf}{\sc} 

\theoremstyle{plain}
\newtheorem{thm}[pulse]{\thf Theorem}

\newtheorem{prop}[pulse]{\thf Proposition}
\newtheorem{lem}[pulse]{\thf Lemma}
\newtheorem{cor}[pulse]{\thf Corollary}
\theoremstyle{definition}
\newtheorem{dfn}[pulse]{\thf Definition}

\newtheorem{eg}[pulse]{\thf Example}
\theoremstyle{remark}
\newtheorem{rem}[pulse]{\thf Remark}

\newenvironment{romnum}{%
\begin{enumerate}
\renewcommand{\theenumi}{{\rm(\roman{enumi})}}
\renewcommand{\labelenumi}{\theenumi}
}{\end{enumerate}\ignorespacesafterend}

\newcommand{\defeq}{:=}
\newcommand{\dt}[1]{\textit{#1}\/}  
\newcommand{\st}{\mathbin{\colon}}

\newcommand{\zer}{\boldsymbol{0}}

\newcommand{\Hom}{\operatorname{Hom}}
\newcommand{\defect}{\operatorname{def}}
\newcommand{\dist}{\operatorname{dist}}

\newcommand{\abs}[1]{{\left\vert#1\right\vert}}

\newcommand{\set}[1]{\left\{#1\right\}}

\newcommand{\join}{\operatorname{join}}

\newcommand{\TWO}{\underline{\sf{2}}}

\newcommand{\Sch}{\operatorname{Sch}}
\newcommand{\Shat}{\widehat{S}} 

\newcommand{\filters}[1]{\operatorname{Filt}(#1)} 
\newcommand{\filgen}[1]{\operatorname{fil}(#1)}  

\newcommand{\br}{\operatorname{br}}  

\newcommand{\FIN}{\operatorname{\mathcal{P}}\nolimits^{\rm fin}} 

\newcommand{\FBP}{\mathop{\rm fbp}\nolimits} 
\newcommand{\stable}[1]{{$\FBP_{#1}$-stable}}
\newcommand{\red}{\operatorname{{\sf r}}}

\newcommand{\veps}{\varepsilon}
\newcommand{\gm}{\gamma}
\newcommand{\lm}{\lambda}
\newcommand{\om}{\omega}

\newcommand{\Om}{\Omega}

\newcommand{\ssfont}[1]{{\mathsf{#1}}}
\newcommand{\ssa}{\ssfont{a}}
\newcommand{\ssb}{\ssfont{b}}

\newcommand{\ssx}{\ssfont{x}}
\newcommand{\ssy}{\ssfont{y}}
\newcommand{\ssz}{\ssfont{z}}

\newcommand{\cE}{{\mathcal E}}
\newcommand{\cF}{{\mathcal F}}

\newcommand{\cP}{{\mathcal P}}

\newcommand{\cS}{{\mathcal S}}

\newcommand{\Cplx}{{\mathbb C}}
\newcommand{\Nat}{{\mathbb N}}

\newcommand{\Torus}{{\mathbb T}}

\newcommand{\sscap}{\wedge} 

\newcommand{\setcap}{\cap} 
\newcommand{\setcup}{\cup} 
\newcommand{\bigsetcup}{\bigcup} 

\renewcommand{\sscap}{\cap}  
\numberwithin{equation}{section}

\usepackage[colorlinks]{hyperref}
\usepackage{sectsty}
\allsectionsfont{\sffamily}

\begin{document}
\maketitle

\begin{abstract}
We continue the study of the AMNM property for weighted semilattices that was initiated in \cite{YC_wtsl-amnm}. We reformulate this in terms of stability of filters with respect to a given weight function, and then provide a combinatorial condition which is necessary and sufficient for this ``filter stability'' property to hold. Examples are given to show that this new condition allows for easier and unified proofs of some results in \cite{YC_wtsl-amnm}, and furthermore allows us to verify the AMNM property in situations not covered by the results of that paper. As a final application, we show that for a large class of semilattices, arising naturally as union-closed set systems, one can always construct weights for which the AMNM property fails.


\end{abstract}

\begin{section}{Introduction}
\subsection{Background context}
This paper is a sequel to \cite{YC_wtsl-amnm}, in which the first author studied approximately multiplicative functionals on weighted convolution algebras of semilattices, and the problem of whether all such functionals arise as small-norm perturbations of multiplicative ones.

Banach algebras with such a ``stability'' property are said to be \dt{AMNM} \cite{BEJ_AMNM1} or \dt{f-stable} \cite{Jarosz_AMNM}.
When $S$ is a semigroup and $\om$ is a (submultiplicative) weight on $S$, the AMNM condition on the Banach algebra $\ell^1_\om(S)$ is equivalent to a kind of ``Hyers--Ulam''-type stability for $\Cplx$-valued functions on $S$, where the stability is measured in terms of the weight~$\om$. Hence, one can define the notion of AMNM for the pair $(S,\om)$ directly without relying on the Banach-algebraic concept (see Definition \ref{d:semigroup AMNM} and Remark \ref{r:AMNM for semigroup or algebra}). Indeed, the present paper has been written so that no Banach algebra theory is required.

In the case where $S$ is a semilattice, a sufficient condition for $(S,\om)$ to be AMNM was given in \cite[Theorem 3.14]{YC_wtsl-amnm}, which covered many natural examples, such as semilattices with finite breadth. However, it was pointed out by the third author in an unpublished communication, that this condition fails to hold for some straightforward examples, such as Example \ref{eg:easy example} below, and so it was clear that the earlier paper left significant room for improvement. For example, is the weighted semilattice in Example \ref{eg:easy example} AMNM?

\subsection{New work}
In this paper we sharpen some of the techniques introduced in \cite{YC_wtsl-amnm}, and put them in a more systematic framework, to obtain a necessary and sufficient criterion for a weighted semilattice $(S,\om)$ to be AMNM.

As an intermediate step, we observe that it suffices to restrict attention to functions on $S$ with values in $\{0,1\}$, rather than those with codomain $\Cplx$.
This was already noted in \cite{YC_wtsl-amnm} but we take the opportunity to give some extra details and make some minor corrections. These and other preliminary results are given in Section~\ref{s:prelim}, which also has a brief discussion of characters on semilattices and the corresponding filters.

The discretization procedure mentioned above allows us to change perspective from ``character stability'' to ``filter stability'' (with respect to the given weight $\om$).
Our criterion is then an intrinsic characterization of the filter stability property, in terms of a combinatorial property of the pair $(S,\om)$ that we call \dt{propagation}. The precise statement is given in Theorem~\ref{t:equivalence}, and most of Section~\ref{s:criterion} is devoted to setting up the necessary framework.

We highlight two easy applications of Theorem \ref{t:equivalence}:
\begin{enumerate}
\renewcommand{\theenumi}{(\arabic{enumi})}
\renewcommand{\labelenumi}{(\arabic{enumi})}
\item
It is easy to check that the weighted semilattice described in Example~\ref{eg:easy example} has ``propagation at all levels''. Therefore, by  Theorem~\ref{t:equivalence}, this weighted semilattice satisfies filter stability, and hence is AMNM. Details are given in Section~\ref{ss:AMNM but not flighty}, which also explains why the results of \cite{YC_wtsl-amnm} are insufficient here.
\item\label{li:fin breadth always AMNM}
It is also easy to check that if $S$ is a semilattice with \dt{finite breadth}, then for \emph{any choice of weight $\om$} we have propagation at all levels. Hence we can apply Theorem \ref{t:equivalence} once again, to deduce that $(S,\om)$ is AMNM. This was already observed in \cite{YC_wtsl-amnm}, but Theorem \ref{t:equivalence} provides a unified framework to see why this result holds.
\end{enumerate}

A third application, which requires substantially more work, forms the final main result of our paper (Theorem \ref{t:subslatt of Pfin}); it addresses a question raised in the closing remarks of \cite{YC_wtsl-amnm}. Given the result stated above in \ref{li:fin breadth always AMNM}, and given that there is a weighted semilattice $(T,\om)$ which is not AMNM (\cite[Theorem 3.4]{YC_wtsl-amnm}), it is natural to ask if \emph{every} $S$ with infinite breadth can be equipped with \emph{some} weight $\om$ such that $(S,\om)$ is not AMNM.
While we are unable to resolve this question here, we use the ``propagation'' criterion to obtain a positive answer for a large class of semilattices: namely, those which embed homomorphically into $(\FIN(\Omega),\cup)$ for some set $\Omega$ (the notation is explained at the start of Section~\ref{s:prelim}).

\subsection{Future work}
We would like to extend the construction used to prove Theorem~\ref{t:subslatt of Pfin} to cover all semilattices, not just those $S$ that embed into $(\FIN(\Om),\cup)$. This runs into some serious technical obstacles. Most notably, our construction relies on being able to locate finite free subsemilattices of $S$ that fit togther in a well-behaved way; this is possible here because finite free subsemilattices of $(\FIN(\Om),\cup)$ only live on a ``finite part'' of $\Om$, and hence can be inductively removed without interfering too much with the embedded copy of~$S$.

To move from $\FIN(\Om)$ to the full powerset $\cP(\Om)$, it seems that a deeper study is required of union-closed set systems with infinite breadth,
and their possible substructures.
We intend to pursue this in greater depth in forthcoming work.

\end{section}

\begin{section}{Background and preliminaries}
\label{s:prelim}

\subsection{Notational conventions}
If $A$ and $B$ are non-empty sets then $A^B$ denotes the set of functions $B\to A$; $\cP(B)$ denotes the powerset of $B$; and $\FIN(B)$ denotes the set of all finite subsets of~$B$.

Some of our later examples are subsets of $\FIN(\Om)$ for some non-empty set $\Om$. When working with those examples, we shall use upper-case calligraphic letters to denote subsets of $\FIN(\Om)$, while using letters such as $\ssa$ or $E$ for subsets of $\Om$; elements of $\Om$ will be denoted by lower-case Greek letters.%

If $S$ and $T$ are semigroups then the set of homomorphisms from $S$ to $T$ is denoted by $\Hom(S,T)$. The zero function on a given set will usually be denoted by $\zer$, the domain being clear from context.

\begin{subsection}{Semilattices, characters and filters}
We recall some standard definitions, for readers who are not semigroup theorists.

A \dt{semilattice} is a commutative semigroup $S$ satisfying $x^2=x$ for all $x\in S$. Such an $S$ has a standard and canonical partial order: if $x,y\in S$ we write $x\preceq y$ whenever $xy=x$.
It is sometimes useful to read $x\preceq y$ as: ``$x$ is a multiple of $y$'' or ``$y$ is a factor of~$x$''. In this language, $xy$ is the ``largest common multiple'' of $x$ and $y$.

With respect to $\preceq$, $xy$ is the \dt{meet} (or~\dt{greatest lower bound}) of $x$ and~$y$; this gives an alternative, order-theoretic definition of a semilattice, as a poset in which every pair of elements has a meet.
Both the semigroup definition and the order-theoretic definition fit well with the intuition from the next example.

\begin{eg}[The $2$-element semilattice]\label{eg:TWO and its powers}
The set $\{0,1\}$, equipped with usual multiplication, is a semilattice, which we shall denote by $\TWO$. (Note that the partial order defined above satisfies $0\preceq 1$.)

If $X$ is a non-empty index set then $\TWO^X$, equipped with co-ordinatewise multiplication, is also a semilattice. Each $\psi\in\TWO^X$ can be identified with $\psi^{-1}(1)\subseteq X$, and this defines an isomorphism of semilattices $\TWO^X \cong (\cP(X),\cap)$.
\end{eg}

Semilattices of the form $(\cP(X),\cap)$ are universal in a certain sense. Given a semi\-lattice~$S$, for each $x\in S$ let
\begin{equation}
E_x=\{t\in S \colon t\preceq x\}.
\end{equation}
The definition of $\preceq$ implies that $E_x\cap E_y= E_{xy}$ for all $x,y\in S$, and so the function $\Sch: S\to (\cP(S),\cap)$ defined by $x\mapsto E_x$ is a semigroup homomorphism. Since $x\in E_x$ for all $x\in S$, it is also easily verified that $\Sch$ is injective. Thus every semilattice can be represented faithfully as a \emph{set system} which is closed under binary intersections --- or, by taking complements, as a set system which is closed under binary unions.
The map $\Sch$ has appeared in the literature under various names, such as the \dt{Cayley embedding} or the \dt{Sch\"utzenberger representation}
\cite{GHS_slatt}.

We now turn to characters. 
Let $S$ be a commutative semigroup, and regard $\Cplx$ as a semigroup with respect to multiplication. Following the terminology of \cite{HZ_TAMS56}, a \dt{semicharacter on $S$} is a semigroup homomorphism $S\to\Cplx$ which is bounded (equivalently, which takes values in the closed unit disc) and which is not identically zero. The set of semicharacters on $S$ is denoted by~$\Shat$.
If $S$ is a semilattice, every semigroup homomorphism $S\to\Cplx$ takes values in $\{0,1\}$, so using the notation of Example~\ref{eg:TWO and its powers} we have
\begin{equation}
\Shat\cup\{\zer\}=\Hom(S,\TWO) \subseteq \TWO^S\;.
\end{equation}
It turns out that $\Shat$ is a sub-semilattice of $\TWO^S$.
This is reminiscent of the situation for abelian groups: non-zero semigroup homomorphisms from an abelian group $G$ to $\Cplx$ take values in $\Torus$, and it turns out that the set of characters $\widehat{G}$ is a subgroup of $\Torus^G$.

\begin{rem}[Conflicting terminology]
In (semi)lattice theory it is more common to call a non-zero homomorphism $S\to\TWO$ a \dt{character} of the semilattice~$S$. We will adopt this convention, which fits \cite{HMS_slatt-dual-preview,HMS_slatt-dual}; the remarks above show that the concepts of semicharacter and character coincide for semilattices. We shall reserve ``semicharacter'' for general commutative semigroups.
\end{rem}

\begin{rem}[Duality theory for semilattices]
\label{r:duality}
There is a duality theory for semilattices in which $\TWO$ acts as the dualizing object, analogous to Pontrjagin duality for (locally compact) abelian groups where $\Torus$ is the dualizing object. See \cite{HMS_slatt-dual} for further details, with an emphasis on category-theoretic aspects; an accessible pr\'ecis, without proofs, is provided in \cite{HMS_slatt-dual-preview}.
Here, we merely remark that there is a natural embedding $g:S\to \TWO^{\Shat}$ which sends $x\in S$ to $\widehat{x}: \phi \mapsto\phi(x)$; this is essentially the same as the Gelfand representation for the commutative Banach algebra $\ell^1(S)$.
\end{rem}

As noted above, $\{0,1\}$-valued functions on $S$ correspond to subsets of $S$. It is standard knowledge in semilattice theory that, under this correspondence, characters on $S$ (viewed as elements of $\Hom(S,\TWO)$) correspond to filters in~$S$. This perspective was already exploited in \cite{YC_wtsl-amnm}, and will be pursued more systematically in this paper. We therefore review some of the details for the reader's convenience.

Traditionally, \dt{a filter} in a semilattice $S$ is defined to be a non-empty subset $F$ which is a subsemigroup (so $x,y\in F\implies xy\in F$) and is ``upwards-closed'' with respect to $\preceq$ (so $x\in S$, $y\in F$ and $xy=y$ $\implies x\in F$). The second condition is equivalent to $xy\in F \implies x,y\in F$ --- we leave the proof to the reader --- and so for this paper we adopt the following equivalent definition, which makes the link with characters clearer.

\begin{dfn}[Filters in semilattices]
  \label{d:filter}
  Let $S$ be a semilattice and let $F\subseteq S$. We say that $F$ is a \dt{filter in $S$} if it is non-empty and satisfies
\[ \forall\; x,y \in S\quad (xy \in F \Longleftrightarrow x,y\in F). \]
The set of all filters in $S$ will be denoted by $\filters{S}$.
\end{dfn}

It follows that for each $X\subseteq S$, $X$ is a filter in $S$ if and only if $1_X\in\Shat$; equivalently,
\begin{equation}\label{eq:character-filter}
\Shat = \{ \psi \in\TWO^S \colon \text{$\psi^{-1}(1)$ is a filter in $S$}\}. 
\end{equation}

We can use duality theory and the link with filters to give an explanation of the Cayley embedding/Sch\"utzenberger representation that was mentioned earlier. Each $y\in S$ defines a filter $F_y=\{t\in S\colon t\succeq y\}$, and so $y\mapsto 1_{F_y}$ defines
an injective function $f:S\to\Shat$.
By restriction, we obtain a homomorphism of semilattices $f^*: \TWO^{\Shat} \to \TWO^S$. Then, composing $f^*$ with
the homomorphism $g: S\to \TWO^{\Shat}$ from Remark \ref{r:duality}
yields
$f^*g: S\to \TWO^S\equiv (\cP(S),\cap)$.
Checking through the definitions, we find that for each $x\in S$
\[
f^*g(x)= \{ y\in S \colon \widehat{x}(1_{F_y})=1\} = \{ y\in S \colon x\in F_y\} = \{ y\in S \colon y\preceq x\}.
\]
Thus $f^*g$ is exactly the embedding $\Sch: S\to (\cP(S),\cap)$ from before.
\end{subsection}

\begin{subsection}{AMNM for weighted semilattices}
In this section, we set up the notion of AMNM for weighted semilattices in a way that does not require definitions or results from the theory of Banach algebras. We also take the opportunity to make precise some arguments from \cite{YC_wtsl-amnm} that were unclear or incomplete as stated there:
see Lemma \ref{l:auto omega-bdd} and Proposition~\ref{p:discretize} below.

We start in a more general setting.
Let $S$ be a semigroup.
A \dt{weight function on $S$} is a function $\om:S\to(0,\infty)$ which is 
\dt{submultiplicative}; that is, $\om(xy)\leq\om(x)\om(y)$ for all $x,y\in S$.
We shall often refer to the pair $(S,\om)$ as a \dt{weighted semigroup.}

Given a weighted semigroup $(S,\om)$ and $\phi,\psi\in\Cplx^S$, we define:
\begin{equation}
\begin{aligned}
\defect_\om(\psi) & = \sup\{ |\psi(xy)-\psi(x)\psi(y)|\om(x)^{-1}\om(y)^{-1} \colon x,y\in S\} \\
d_\om(\psi,\phi) & = \sup_{x\in S} |\psi(x)-\phi(x)|\om(x)^{-1} \\
\dist_\om(\psi) & = \inf\{ d_\om(\psi,\phi) \colon \phi\in\Hom(S,\Cplx) \}
\end{aligned}
\end{equation}
where in the last definition we view $\Cplx$ as a \emph{multiplicative} semigroup.
Note that in all three definitions we allow the value $+\infty$.

$\defect_\om(\psi)$ and $\dist_\om(\psi)$ give two different ways to quantify the failure of $\psi$ to be a homomorphism;
the first has a ``local'' flavour and the second has a ``global'' flavour.
We shall sometimes refer to $\defect_\om(\psi)$ as the \dt{$\om$-defect} of~$\psi$.

As in \cite{YC_wtsl-amnm}, we say that $\psi\in\Cplx^S$ is \dt{$\om$-bounded} if
$\sup_{x\in S} |\psi(x)|\om(x)^{-1}<\infty$.
The following lemma is not needed for the eventual application to semilattices, but we include it for sake of completeness and for possible use in future work.

\begin{lem}[Improving $\om$-boundedness]
\label{l:improve omega-bdd}
Let $S$ be a semigroup. If $\psi\in\Cplx^S$ is $\om$-bounded, then
$\sup_{x\in S}|\psi(x)|\om(x)^{-1} \leq 1 + \defect_\om(\psi)$.
\end{lem}

This follows from a technique known in the Banach-algebraic setting \cite[Prop.\ 5.5]{Jarosz_LNM}. To make this paper more self-contained we provide a full proof, adapted from parts of the proof in \cite{Jarosz_LNM}.

\begin{proof}[Proof (Jarosz).]
Let $K= \sup_{x\in S} |\psi(x)|\om(x)^{-1}$, and let $s\in S$. Since
\[ \frac{|\psi(s)^2 - \psi(s^2)|}{\om(s)^2} \leq \defect_\om(\psi) \]
we have
\[ \frac{|\psi(s)|^2}{\om(s)^2} -\defect_\om(\psi) \leq \frac{|\psi(s^2)|}{\om(s)^2} \leq \frac{|\psi(s^2)|}{\om(s^2)} \leq K.\]
Taking the supremum over all $s$ yields $K^2-\defect_\om(\psi)\leq K$. Since $t-1\leq t^2-t$ for all $t\geq 0$, we obtain $K-1\leq K^2-K \leq \defect_\om(\psi)$ as required.
\end{proof}

The lemma may fail if we drop the assumption that $\psi$ is $\om$-bounded.
For instance, take $S=\Nat$ with usual addition, take $\om\equiv 1$, and let $\psi(n)=2^n$. Then $\psi$ is clearly not $\om$-bounded, yet $\defect_\om(\psi)=0$.

The following terminology is modelled on the corresponding terminology for Banach algebras, as found in \cite{Jarosz_AMNM} or \cite{BEJ_AMNM1}.

\begin{dfn}[AMNM/f-stability for weighted semigroups]
\label{d:semigroup AMNM}
Let $(S,\om)$ be a weighted semigroup. We say that $(S,\om)$ has the \dt{AMNM property}, or is \dt{AMNM}, or is \dt{f-stable}, if the following holds: for all $\veps>0$ there exists $\delta>0$ such that every $\om$-bounded $\psi\in\Cplx^S$ with $\defect_\om(\psi) < \delta$ satisfies $\dist_\om(\psi) <\veps$.
\end{dfn}

\begin{rem}\label{r:AMNM for semigroup or algebra}
Given a weighted semigroup $(S,\om)$, we can form the associated weighted convolution algebra $\ell^1_\om(S)$. Then $(S,\om)$ is AMNM if and only if $\ell^1_\om(S)$ has the Banach-algebraic AMNM property of \cite{BEJ_AMNM1}. The proof of this equivalence is merely an exercise in translating the definitions, together with some standard properties of the $\ell^1$-norm; details can be found in \cite[Section~2.2]{YC_wtsl-amnm}.
\end{rem}

We now specialize to the setting where $S$ is a semilattice. Note first that if $\om$ is a weight function on a semilattice, we have $\om(x)^2\geq\om(x^2)=\om(x)$ for all $x\in S$, which forces $\om(x)\geq 1$ for all $x\in S$.

Another special feature is that we can improve on Lemma \ref{l:improve omega-bdd}. The following lemma plugs a small gap in the statement/proof of \cite[Corollary 3.7]{YC_wtsl-amnm}.

\begin{lem}\label{l:auto omega-bdd}
Let $(S,\om)$ be a weighted semilattice and let $\psi\in\Cplx^S$. Then\newline
$\sup_{x\in S} |\psi(x)|\om(x)^{-1} \leq 1+\defect_\om(\psi)$.
In particular, if $\defect_\om(\psi)<\infty$ then $\psi$ is $\om$-bounded.
\end{lem}

\begin{proof}
Without loss of generality, assume $\defect_\om(\psi)<\infty$. Let $s\in S$. Arguing similarly to the proof of Lemma \ref{l:improve omega-bdd}, we obtain
\[
|\psi(s)|\om(s)^{-1} = |\psi(s^2)|\om(s)^{-1}
\geq |\psi(s^2)|\om(s)^{-2}
\geq |\psi(s)^2| \om(s)^{-2} - \defect_\om(\psi)
\]
and hence (again following the proof of the earlier lemma)
\[ \frac{|\psi(s)|}{\om(s)} -1 \leq \left(\frac{|\psi(s)|}{\om(s)}\right)^2 - \frac{|\psi(s)|}{\om(s)} \leq \defect_\om(\psi) \]
as required.
\end{proof}

For $X,Y\subseteq S$, we shall abbreviate $\defect_\om(1_X)$, $d_\om(1_X,1_Y)$ and $\dist_\om(1_X)$ to $\defect_\om(X)$, $d_\om(X,Y)$ and $\dist_\om(1_X)$ respectively. By the remarks following Definition \ref{d:filter},
\[
\dist_\om(X) = \inf\{ d_\om(X,F) \colon F\in \filters{S}\cup \emptyset
\}.
\]

The following result is essentially the same as \cite[Corollary 3.7]{YC_wtsl-amnm}, phrased in different language. However, the proof there was left to the reader, and as mentioned above there was a missing step glossed over. 
 We take this opportunity to provide a slightly weaker but more precise statement, with a complete proof.

\begin{prop}[Discretizing the AMNM problem for weighted semilattices]\label{p:discretize}
Let $(S,\om)$ be a weighted semilattice. The following are equivalent:
\begin{romnum}
\item\label{li:original}
$(S,\om)$ is AMNM;
\item \label{li:epsilon-delta}
$\forall\ \veps>0\ \exists\ \delta>0$ such that every $G\subseteq S$ with $\defect_\om(G)< \delta$ satisfies $\dist_\om(G)<\veps$.
\item\label{li:sequences}
every sequence $(X_n)$ in $\cP(S)$ with $\defect_\om(X_n)\to 0$ satisfies $\dist_\omega(X_n)\to 0$. 
\end{romnum}
\end{prop}

\begin{proof}
The equivalence of \ref{li:epsilon-delta} and \ref{li:sequences} is routine.
The implication \ref{li:original}$\implies$\ref{li:epsilon-delta} follows directly from the definitions
and the fact that $1_G$ is $\om$-bounded for all $G\subseteq S$ (since $\om\geq1$).
%
%

Finally, suppose \ref{li:epsilon-delta} holds. Let $\veps>0$, and let $\delta=\delta(\veps)$ be as provided in~\ref{li:epsilon-delta}. Let $\delta_1>0$, to be determined later, and let $\psi\in\Cplx^S$ be $\om$-bounded with $\defect_\om(\psi) <\delta_1$. By \cite[Lemma~3.6]{YC_wtsl-amnm} and the calculations which follow it, there is a set $G\subseteq S$ such that $d_\om(\psi,1_G) < {\delta_1}^{1/2}$ with $\defect_\om(G) \leq 3{\delta_1}^{1/2}+2\delta_1$. Therefore, provided we originally chose $\delta_1$ small enough that $\delta_1< \veps^2$ and $3{\delta_1}^{1/2}+2\delta_1 < \delta(\veps)$, condition~\ref{li:epsilon-delta}
implies there exists $\phi\in\Shat\cup\{\zer\}$ such that
\[ d_\om(\psi,\phi) \leq d_\om(\psi, 1_G) + d_\om(1_G,\phi) <2\veps. \]
Thus \ref{li:original} holds.
\end{proof}

Conditions \ref{li:epsilon-delta} and \ref{li:sequences} can be viewed as a kind of ``filter stability'' property relative to the weight function $\om$, and for the rest of this paper we shall work exclusively with the perspective of subsets of $S$ rather than functions on $S$.

\end{subsection}
\end{section}


\begin{section}{Characterizing stability of filters}
\label{s:criterion}

\subsection{Initial remarks}
In this section we shall derive an intrinsic, combinatorial characterization of when a given $(S,\om)$ has ``stable filters'' in the sense of Proposition \ref{p:discretize}\ref{li:epsilon-delta}.
The precise statement of our characterization is given
as Theorem~\ref{t:equivalence}, and we shall build up to it in stages.

It will be convenient for later examples if we switch to working with \dt{log-weights}, by which we mean functions $\lm:S \to [0,\infty)$ that satisfy $\lm(xy)\leq \lm(x)+\lm(y)$ for all $x,y\in S$. Given such a $\lm$ and $L\geq 0$, we define
$W_L(S,\lm)=\{x\in S \colon \lm(x)\leq L\}$.
When there is no danger of confusion we abbreviate this to $W_L$.

\begin{eg}\label{eg:easy example}
Let $\Omega$ be a non-empty set. $\FIN(\Omega)$ is a semilattice with respect to binary union.
The function $\ssx\mapsto |\ssx|$ is a log-weight on $(\FIN(\Omega),\cup)$, and
\[
W_1 = \{\emptyset\}\cup \set{ \{\gm\} \colon \gm\in \Omega}.
\]
If we fix $r\in (1,\infty)$, the function $\om_r: \ssx\mapsto r^{|\ssx|}$ is a weight on $(\FIN(\Omega),\cup)$.
\end{eg}

By varying the ``level'' $L$ we obtain a useful filtration of~$S$.
Intuitively, we think of two given sets as ``agreeing at level~$L$'' when they have the same intersection with $W_L$.

\subsection{Characterizing sets on which $\dist_\om$ is small.}

As a first step towards Theorem~\ref{t:equivalence}, we show that given two subsets $E_1,E_2\subseteq S$, $d_\om(E_1,E_2)$ is small precisely when $E_1$ and $E_2$ ``agree at level~$L$'' for some large value of~$L$.

\begin{lem}\label{l:levels}
Let $E_1,E_2\subseteq S$ and let $\lm=\log\om$.
\begin{romnum}
\item
Let $\veps>0$. If $d_\om(E_1,E_2)<\veps$, then $E_1\cap W_L=E_2\cap W_L$ for every $L \leq\log(\veps^{-1})$.
\item
Let $L\geq 0$. If $E_1\cap W_L=E_2\cap W_L$ then $d_\om(E_1,E_2)\leq e^{-L}$.
\end{romnum}
\end{lem}

\begin{proof}
In both parts, the key point is that $\abs{1_{E_1}(x)-1_{E_2}(x)}$ takes values $0$ or $1$.
\begin{romnum}
\item
Let $L \leq\log(\veps^{-1})$ and suppose $d_\om(E_1,E_2)<\veps$. Then $\abs{1_{E_1}(x)-1_{E_2}(x)} < \veps\om(x)$ for all $x\in S$.
For each $x\in W_L$, $\veps\om(x) \leq 1$ and so $1_{E_1}(x)=1_{E_2}(x)$; thus $E_1\cap W_L=E_2\cap W_L$.
\item
Suppose
$E_1\cap W_L=E_2\cap W_L$. Then $1_{E_1}(x)=1_{E_2}(x)$
whenever $\om(x)\leq e^L$. Hence $\abs{1_{E_1}(x)-1_{E_2}(x)} \leq\min(1, e^{-L}\om(x))$ for all $x$, which implies $d_\om(E_1,E_2)\leq e^{-L}$.
\end{romnum}
\vskip-2.0em
\end{proof}

If $E\subseteq S$, let $\filgen{E}$ denote the \dt{filter-or-empty-set generated by~$E$}, i.e.~the intersection of all $X\subseteq S$ such that $X\in\filters{S}\cup\{\emptyset\}$ and $X\supseteq E$. Note that $\filgen{E}=\emptyset$ if and only if $E=\emptyset$.

\begin{rem}\label{r:best guess}
Let $E\subseteq S$ and let $L\geq 0$. Suppose there is some $F\in\filters{S}\cup\{\emptyset\}$ such that $E\cap W_L=F\cap W_L$.
Since $F$ contains $E\cap W_L$ it contains $\filgen{E\cap W_L}$. Hence
\[
E\cap W_L = F\cap W_L \supseteq \filgen{E\cap W_L} \cap W_L \supseteq E\cap W_L
\]
which forces $E\cap W_L = \filgen{E\cap W_L}\cap W_L$.

Summarizing: if $E$ agrees with \emph{some} filter-or-empty-set at level $L$, it agrees with the \emph{particular} filter-or-empty-set $\filgen{E\cap W_L}$ at level~$L$.
Hence, in view of Lemma~\ref{l:levels}, $\dist_\om(E)$ is small if and only if $E$ agrees with $\filgen{E\cap W_L}$ at level $L$ for some large $L$.
\end{rem}

\subsection{Characterizing sets with small $\om$-defect}

Let $X\subseteq S$ be non-empty. Note that if $x,y\in X$ and $z\succeq xy$ (that is, $z$ is a factor of $xy$) then $z\in\filgen{X}$. Moreover, every $z\in \filgen{X}$ satisfies $z\succeq x_1\dotsb x_k$ for some $x_1,\dots, x_k\in X$.

\begin{dfn}[$\FBP_C$-stability]
Let $C\geq 0$.
For $X\subseteq S$ we define
\begin{equation}
\FBP_C(X)\defeq \set{ z\in W_C \colon \hbox{there exist $x,y \in X\cap W_C$ such that $z\succeq xy$} } .
\end{equation}
($\FBP$ stands for ``factors of binary products''.)
Note that $\FBP_C(\emptyset)=\emptyset$. 
We always have $\FBP_C(X)\supseteq X\cap W_C$; if equality holds, we say that $X$ is \dt{\stable{C}}.
\end{dfn}

$X$ is \stable{C} if and only if $\FBP_C(X)\subseteq X$: the ``only if'' direction is trivial, and the ``if'' direction holds because $\FBP_C(X)\subseteq W_C$ for any $X\subseteq S$.

\begin{lem}\label{l:stable-and-AM}
Let $(S,\om)$ be a weighted semilattice, and put $\lm=\log\om$.
Let $X\subseteq S$.
\begin{romnum}
\item\label{li:stable-to-AM}
Let $C\geq 0$. If $X$ is \stable{C}, then $\defect_\om(X) \leq e^{-C}$.
\item\label{li:AM-to-stable}
Let $\delta>0$. If $\defect_\om(X)< \delta$ then $X$ is \stable{C} for any $C\leq \log(\delta^{-1/3})$.
\end{romnum}
\end{lem}

\begin{proof}
\begin{romnum}
\item
Assume $X$ is \stable{C}, and let $x,y\in S$. We consider two cases.

Case A: $\lm(x)+\lm(y)\geq C$. In this case
\[
1 \leq e^{-C} e^{\lm(x)+\lm(y)} = e^{-C} \om(x)\om(y),
\]
and so  $\abs{1_X(xy)-1_X(x)1_X(y)}\leq e^{-C} \om(x)\om(y)$.

 Case B: $\lm(x)+\lm(y)\leq C$. In this case
 $x$, $y$ and $xy$ all belong to $W_C$. Since $X$ is \stable{C}:
  \begin{itemize}
\item[--] if $x$ and $y$ are in $X$, then
 $xy\in\FBP_C(\{x,y\})\subseteq \FBP_C(X) \subseteq X$;
\item[--] if $xy$ lies in $X$, then both $x$ and $y$
 belong to $\FBP_C(\{xy\})\subseteq \FBP_C(X) \subseteq X$.
\end{itemize}
Thus in this case, $1_{X}(xy)=1_X(x)1_X(y)$.

Putting Case A and Case B together, we see that $\defect_\om(X)\leq e^{-C}$.

\item
If $X\cap W_C=\emptyset$ then we are done. So suppose $X\cap W_C\neq \emptyset$.
Let $x,y\in X\cap W_C$,
and let $z\in W_C$ with $z\succeq xy$.
Then
\begin{equation}\label{eq:miscell}
\begin{aligned}
\abs{1_X(xy)- 1_X(x)1_X(y) } & \leq \defect_\om(X) \om(x)\om(y) & \leq \defect_\om(X) e^{2C} , \\
\abs{1_X(xy)- 1_X(xy)1_X(z) } & \leq \defect_\om(X)\om(xy)\om(z) &\leq \defect_\om(X) e^{3C} ,
\end{aligned}
\tag{$*$}
\end{equation}
and both are strictly less than $1$ if
$\defect_\om(X) < \delta$ and
$C\leq \log(\delta^{-1/3})$.

Since $1_X(x)=1=1_X(y)$, the first formula in \eqref{eq:miscell} implies $1_X(xy)=1$; feeding this into the second formula yields $z\in X$. Therefore $\FBP_C(X\cap W_C)\subseteq X\cap W_C$, and the converse inclusion is trivial.
\end{romnum}
\vskip-1.0em
\end{proof}

Summarizing: a subset of $S$ has small $\om$-defect precisely when it is \stable{C} for some large value of~$C$.

\subsection{$\FBP_C$-stability and propagation}
An obvious way to obtain \stable{C} sets is by iteration.
Given $E\subseteq S$, put $\FBP_C^0(E)=E\cap W_C$ and for $k\geq 1$
recursively define $\FBP_C^k(E)=\FBP_C(\FBP_C^{k-1}(E))$. Then
\[ 
	E\cap W_C \subseteq \FBP_C(E) \subseteq \FBP_C^2(E) \subseteq \FBP_C^3(E)\subseteq \dots
\]
Define $\FBP_C^\infty(E)$ to be the inductive limit $\bigcup_{k\geq 1} \FBP_C(E)$.
By induction,
\begin{equation}\label{eq:inclusion}
E\cap W_C \subseteq \FBP_C^\infty(E) \subseteq \filgen{E}\cap W_C \quad\hbox{for all $E\subseteq S$.}
\end{equation}%

Even when $S$ is finite, the second inclusion in \eqref{eq:inclusion} can be proper, as we will see in Example~\ref{eg:prototype}.
On the other hand, if $E$ is non-empty and $z\in \filgen{E}$, there always exists some $C\geq 0$, possibly depending on $z$, such that $z\in \FBP_C^\infty(E)$.
(For instance, if $x_1,\dots, x_k\in E$ and $z\succeq x_1\dotsb x_k$ then $C=\max\{\sum_{i=1}^k\lm(x_i), \lm(z)\}$ suffices.)
This leads naturally to the following definition.

\begin{dfn}[Propagation]
\label{d:propagation}
For $z\in\filgen{E}$, let
\begin{equation}
V_E(z) = \inf \set{ C\geq 0 \colon z\in \FBP_C^{\infty}(E)}.
\end{equation}
Given $L\geq 0$, we say that $(S,\lm)$ \dt{propagates at level $L$}, or \dt{has $L$-propagation}, if
\begin{equation}
\label{eq:sup sup finite}
\sup_{\emptyset\neq E \subseteq W_L} \sup_{z\in \filgen{E}\cap W_L}V _E(z) < \infty\;.
\end{equation}
It is convenient to  set $V_E(z)\defeq+\infty$ whenever $z\notin\filgen{E}$.
\end{dfn}

We make some observations for future reference:
\begin{itemize}
\item
For every $z\in S$ we have $V_E(z)\geq \lm(z)$
and $V_E(z)\geq \inf_{x\in E} \lm(x)$.
(The first inequality holds since $\FBP_C^\infty(E)\subseteq W_C$. The second one holds because if $C< \inf_{x\in E} \lm(x)$ then $\FBP_C(E)=\emptyset$, preventing $z\in\FBP_C^\infty(E)$.)
\item
If $(S,\lm)$ propagates at a level $L$, then it also does so at every lower level.
\item
In the formula defining $L$-propagation, we could restrict $E$ to the finite subsets of $W_L$ without altering the value of the double supremum. 
\end{itemize}

We can now state and prove the promised characterization of ``filter stability''.
\begin{thm}\label{t:equivalence}
Let $(S,\om)$ be a weighted semilattice and let $\lm=\log\om$. The following conditions are equivalent. 
\begin{romnum}
\item\label{li:(i)}
$(S,\om)$ is AMNM.

\item\label{li:(ii)}
$\forall\ \veps>0\ \exists\ \delta>0$ such that every $G\subseteq S$ with $\defect_\om(G)< \delta$ satisfies $\dist_\om(G)<\veps$.

\item\label{li:(iii)}
$\forall L\geq 0\ \exists\ C\geq 0$ such that every $G\subseteq S$ which is \stable{C} satisfies $G\cap W_L = \filgen{G\cap W_L}\cap W_L$.

\item\label{li:(iv)}
$(S,\lm)$ has $L$-propagation for all $L\geq 0$.
\end{romnum}

\end{thm}

\begin{proof}
The equivalence of \ref{li:(i)} and \ref{li:(ii)} was proved in Proposition \ref{p:discretize}.
The equivalence of \ref{li:(ii)} and \ref{li:(iii)} has been demonstrated in the previous two subsections: specifically, one combines Lemma \ref{l:levels}, Remark \ref{r:best guess} and Lemma \ref{l:stable-and-AM} with some standard epsilon-delta book-keeping.

Suppose that condition \ref{li:(iv)} holds. Let $L\geq 0$.
Then we may choose $C\geq L$ such that $V_G(z)\leq C$ for all $G\subseteq W_L$ and all $z\in\filgen{G}\cap W_L$.
Let $E\subseteq S$ be $\FBP_C$-stable and put $G=E\cap W_L$. By our choice of $C$ and the definition of $V_G(z)$,
\[ 
\filgen{E\cap W_L}\cap W_L \subseteq \FBP_C^\infty(E\cap W_L).
\]
But since $E$ is \stable{C} the right-hand side is contained
in~$E$. Hence
$\filgen{E\cap W_L}\cap W_L \subseteq E \cap W_L$,
and the converse inclusion is trivial.
Thus condition \ref{li:(iii)} holds.

Conversely, suppose Condition \ref{li:(iii)} holds.
Let $L\geq 0$ be arbitrary, and choose $C\geq 0$ with the property stated in \ref{li:(iii)}. Let $E\subseteq W_L$.
Put $K=\max(L,C)$. Since $C\leq K$, it follows from the definitions that every \stable{K} set is \stable{C}. In particular, $\FBP_K^\infty(E)$ is \stable{C}, so by the property stated in \ref{li:(iii)},
\[
\FBP_K^\infty(E)\cap W_L = \filgen{ \FBP_K^\infty(E)\cap W_L}\cap W_L.
\]
But since $K\geq L$ and $E\subseteq W_L$, 
the right-hand side contains $\filgen{E}\cap W_L$.
Hence $\filgen{E}\cap W_L \subseteq \FBP_K^\infty(E)$.
This shows that the double supremum in \eqref{eq:sup sup finite} is
$\leq K$,
hence finite; so $(S,\lm)$ has $L$-propagation, and \ref{li:(iv)} holds.
\end{proof}

\subsection{An example where our theorem improves on earlier results}\label{ss:AMNM but not flighty}

In the language of the present paper, \cite[Theorem 3.14]{YC_wtsl-amnm} may be restated as follows. Let $(S,\om)$ be a weighted semilattice, with $\lm=\log\om$, such that
\begin{equation}\label{eq:flighty}
\text{for each $L\geq 0$, $\lm$ is bounded on the subsemilattice generated by $W_L(S,\lm)$.}
\tag{$\diamondsuit$}
\end{equation}
Then $(S,\om)$ is AMNM.

Consider the semilattice $\cS$ and log-weight $\lm$ from Example~\ref{eg:easy example}. That is, $\cS=(\FIN(\Omega),\cup)$ for some nonempty set $\Omega$, and $\lm(\ssx)=|\ssx|$ for $\ssx\in\cS$. Since $\{\gm\}\in W_1$ for each $\gm\in\Om$, the subsemilattice generated by $W_1$ is $\{\ssx\in\cS \colon \ssx\neq\emptyset\}$. Hence the condition in \eqref{eq:flighty} fails for $L=1$, and \cite[Theorem 3.14]{YC_wtsl-amnm} cannot be applied.
On the other hand, we will now show that $(\cS,\lm)$ has $L$-propagation for all $L\geq 0$. Putting $\om=\exp\lm$, it follows from Theorem~\ref{t:equivalence} that $(\cS,\om)$ is AMNM.

Our task will be made slightly simpler using  two preliminary observations:
\begin{itemize}
\item
we only need to verify $L$-propagation for all $L\in\Nat$;
\item
we only need to check the supremum in \eqref{eq:sup sup finite} for \emph{finite} subsets of $W_L$.
\end{itemize}
(For the justification, see the remarks following Definition~\ref{d:propagation}.)

Let $L\in\Nat$, and let $\cE$ be a non-empty finite subset of $W_L$. Put $\ssa = \bigsetcup_{\ssx\in\cE} \ssx$, which is finite and non-empty. A straightforward argument shows that $\filgen{\cE} = \cP(\ssa)$, so that
\[ \filgen{\cE} \sscap W_L = \set{ \ssz\subseteq\ssa \colon |\ssz| \leq L }. \]
If $\ssz\subseteq\ssa$ then there exist $\ssx_1,\dots, \ssx_k\in\cE$ such that $\ssz\subseteq \ssx_1\setcup \dots \setcup \ssx_k$; moreover, if $|\ssz|\leq L$ then we can always choose $k\leq L$. Therefore, $V_{\cE}(\ssz)\leq L^2$ for all $\ssz\in\filgen{\cE}\sscap W_L$. Since this upper bound holds for all finite non-empty $\cE\subseteq W_L$,  $(\cS,\lm)$ propagates at level $L$, just as we claimed.
\end{section}

\begin{section}{Breadth and propagation}
The following is a small variation on the previous example.
\begin{eg}[Free semilattices]
\label{eg:free semilattices}
Let $\Om$ be a non-empty set, and let $\FIN_*(\Om)$ denote the set of finite non-empty subsets of $\Om$.
Equipped with binary union,
$\FIN_*(\Om)$ is a semilattice, which we refer to as the \dt{free semilattice} generated by $\Om$.

This terminology is justified by the following universal property:
for any semilattice $S$ and any function $f:\Om\to S$, there is a unique homomorphism $\widetilde{f}:(\FIN_*(\Om),\cup)\to S$ that extends $f$, defined by $\widetilde{f}(\ssa)=\prod_{\gm\in\ssa} f(\gm)$.
\end{eg}

When $\Om_{00}$ is a finite, non-empty set, we write $\cP_*(\Om_{00})$ instead of $\FIN_*(\Om_{00})$. Semilattices of the form $(\cP_*(\Om_{00}),\cup)$ were used in \cite[Theorem 3.4]{YC_wtsl-amnm} as the building blocks for a weighted semilattice $(T,\om)$ that does not have the AMNM property. The new notion of propagation gives a very natural viewpoint on this construction, with the key details demonstrated in the following example.

\begin{eg}
\label{eg:prototype}
Let $\Omega_{00}$ be a finite set with at least two elements.
 Define
\begin{equation}\label{eq:basic example of bad weight}
\lm(\ssz) \defeq |\ssz| \;\text{ if $\ssz\subsetneq\Omega_{00}$}\quad,\quad \lm(\Omega_{00}) \defeq 0. \qquad(\ssz\in\cP_*(\Omega_{00})).
\end{equation}
Then $\lm$ is a log-weight on the semilattice $(\cP_*(\Om_{00}),\cup)$.
Let $\cE = \set{ \{\om \} \colon \om\in\Omega_{00} }$, and note that 
$W_1 = \cE \cup \{\Omega_{00}\}$.

In particular $\Omega_{00} \in \filgen{\cE}\cap W_1$.
On the other hand, 
$V_\cE(\Omega_{00})\geq \frac{1}{2} |\Omega_{00}|$.
To see this, let $C\geq 0$ be such that $\Omega_{00} \in \FBP_C^\infty(\cE)$.
We have $\Omega_{00}\notin \cE\cap W_C=\FBP_C^0(\cE)$;
let $m\in\Nat$ be minimal such that $\Omega_{00}\in \FBP_C^m(\cE)$. Then there exist $\ssa_1,\ssa_2\in \FBP_C^{m-1}(\cE)$ such that $\Omega_{00}\subseteq \ssa_1\cup\ssa_2$. By minimality of $m$, both $\ssa_1$ and $\ssa_2$ are \emph{proper} subsets of $\Omega_{00}$, and so
\[ |\Omega_{00}| \leq |\ssa_1| + | \ssa_2 | = \lm(\ssa_1)+\lm(\ssa_2) \leq 2C. \]
Hence $2V_E(\Omega_{00}) \geq |\Omega_{00}|$, as claimed.

Intuitively, for this log-weight, the constraint in the $\FBP_C$ operation that we can only multiply elements of log-weight $\leq C$ creates a barrier separating us from~$\Om_{00}$, even though $\Om_{00}$ itself has small log-weight.
\end{eg}

The previous example suggests that to construct log-weights on a given $S$ for which propagation fails, we should look for isomorphic copies of $(\cP_*(\Om_{00}),\cup)$ inside $S$.
If $E\subset S$ is a finite subset, with $\iota_E:E\to S$ being the inclusion map, consider $\widetilde{\iota_E} : \cP_*(E) \to S$, in the notation of Example \ref{eg:free semilattices}.
$\widetilde{\iota_E}$ maps $\cP_*(E)$ onto $\langle E\rangle$, the subsemigroup of~$S$ geneated by~$E$; it is injective if and only if $\langle E\rangle$ has maximal cardinality $|\cP_*(E)|=2^{|E|-1}$.

\begin{dfn}
Let $S$ be a semilattice.
Given a finite, non-empty subset $E\subseteq S$, we say $E$ is \dt{compressible} if there exists a proper subset $E'\subset E$ such that $\prod_{x\in E} x = \prod_{x\in E'} x$; otherwise, we say $E$ is \dt{incompressible}.
\end{dfn}

It is a straightforward exercise to show that $E$ is incompressible if and only if $\widetilde{E}:\cP_*(E)\to S$ is injective (c.f.~Exercise~6(c) of \cite[Section II.5]{Birkhoff_ed3}).

\begin{rem}[Comparison with older terminology]
Our terminology is not entirely standard; the same property is referred to in \cite{LLM_Pac77,Mislove_inproc86} as ``meet irredundant''. However, for examples which arise as subsemilattices of $(\cP(\Om),\setcup)$, as in the next section, the canonical partial order is not given by inclusion but by containment: $\ssa\preceq\ssb \Longleftrightarrow \ssa\supseteq \ssb$; and the ``meet'' of $\ssa$ and $\ssb$ with respect to $\preceq$ is not $\ssa\cap\ssb$, but rather $\ssa\cup\ssb$. In other words, for natural examples we want to consider, the product operation is naturally viewed as a join rather than a meet, and the terminology ``incompressible'' seems more appropriate.
\end{rem}

\begin{dfn}
The \dt{breadth} of a semilattice $S$ is defined to be
\[ \begin{aligned}
\br(S)
& =\sup \set{n\in\Nat \st \text{every subset $E\subseteq S$ with $n+1$ elements is compressible} } \\
& =\sup \set{ |E| \colon \text{$E$ is a finite incompressible subset of $S$}  }
\end{aligned}
. \]
\end{dfn}

 The breadth of a semilattice sheds some light on its structure, and is related to more familiar order-theoretic concepts such as height and width.
 For instance, 
suppose $\br(S)\geq n$. By examining incompressible subsets, one sees that 
 $S$ contains a chain (totally ordered subset) and an antichain (subset in which no two elements are comparable) both of cardinality~$n$;
 see \cite[Section 4.1]{ADHMS_VC2} for further details and some references.
In particular, a semilattice $S$ has breadth~$1$ exactly when the poset $(S,\preceq)$ is totally ordered.

Diverse behaviour occurs even among semilattices of breadth $2$.
For instance, the following example shows that every infinite $k$-ary rooted tree ($k\geq 2$) is a semilattice with breadth $2$ that contains infinite chains and infinite antichains.

\begin{eg}
\label{eg:tree}
Let $k\geq 2$. An infinite $k$-ary rooted tree is an infinite rooted tree in which every vertex has $k$ children. If $x$ and $y$ are vertices in the tree then they have a ``youngest'' common ancestor, which we denote by $x\wedge y$.
 Clearly $\wedge$ is a commutative, associative and idempotent binary operation; so the set of vertices becomes a semilattice $(T,\wedge)$, and the partial order $\preceq$ becomes ``is an ancestor~of''.

There is an infinite path $P\subset T$ obtained by starting at the root and successively taking one of the children; this gives us an infinite chain in $(T,\preceq)$. If for each element of $P$ we then take one of its children that is not in $P$, then the collection obtained is an infinite antichain in $(T,\preceq)$.

Let $x,y,z\in T$, and let $p=x\wedge y\wedge z$. Then either $x\wedge y$ or $y\wedge z$ is equal to $p$: for if not, the set
$\{ p, x\wedge y, y\wedge z, y\}$ would form a cycle of length $\geq 3$ in the tree $T$, which is impossible. Thus every $3$-element subset of $S$ is compressible, and so~$\br(S)\leq 2$. On the other hand, $\br(S)\geq 2$, since $S$ is not totally ordered. 
\end{eg}

\begin{rem}\label{r:finite breadth always AMNM}
Suppose $S$ is a semilattice with finite breadth.
It was originally shown in \cite[Example 3.13 and Theorem 3.14]{YC_wtsl-amnm} that in this case, $(S,\om)$ is AMNM for every weight function $\om$. With the tools of the previous section, we can give an alternative approach to this result. Let $n=\br(S)$, and let $E$ be a finite subset of $S$.

As observed in the remarks before Definition~\ref{d:propagation}: if $x_1,\dots,x_k\in E$ and $z$ is a factor of $x_1\dotsb x_k$, then for any
$C\geq \max\{\sum_{i=1}^k\lm(x_i), \lm(z)\}$
we have $z\in \FBP_C^\infty(E)$. But since $\br(S)=n$, every $z\in \filgen{E}$ is a factor of $x_1\dotsb x_k$ for some $x_1,\dots,x_k\in E$ with $k\leq n$. It follows that whenever $E\subseteq W_L$, $V_E(z)\leq nL$ for all $z\in \filgen{E}\cap W_L$, and therefore $(S,\lm)$ propagates at level $L$ for any $L\geq 0$ and any log-weight $\lm$.  
Now we can apply Theorem~\ref{t:equivalence}.
\end{rem}

Motivated by this result, it is natural to ask the following question,
which was raised implicitly in \cite[Section 6]{YC_wtsl-amnm}.

\paragraph{\sc Question.}
Let $S$ be a semilattice with infinite breadth. Does there exist a log-weight on $S$ for which propagation fails at some level?

\smallskip

At present, we do not have an answer in full generality. However, when $S$ is (isomorphic to) a sub-semilattice of $(\FIN(\Om),\cup)$, the answer is positive.
This is the final main result of our paper.

\begin{thm}\label{t:subslatt of Pfin}
Let $\Omega$ be a non-empty set, and let $\cS\subseteq\FIN(\Om)$ be closed under binary union. If $(\cS,\cup)$ has infinite breadth, then there is a log-weight on $\cS$ which fails $1$-propagation.
\end{thm}


The proof of Theorem~\ref{t:subslatt of Pfin} takes up the rest of this section. From here onwards,
let $\cS$ be a subsemilattice of $(\FIN(\Om),\cup)$ which has infinite breadth.
If~$\cE$ is a finite subset of $\cS$, we write $\join(\cE)$ for $\bigsetcup_{\ssa\in\cE} \ssa \in \cS$.

\begin{lem}\label{l:inductive step}
  Let $\ssa\in\FIN(\Om)$ and let $n \geq 2$. Then there exist $\ssb_1,\dots,\ssb_n\in\cS$ and $\gm_1,\dots,\gm_n\in\Om\setminus \ssa$, such that
  $\gm_j \in \ssb_j \setminus \ssb_k$ whenever $j\neq k$.
\end{lem}

\begin{proof}
Let $M=|\ssa|<\infty$. Since $\cS$ has infinite breadth it has an incompressible subset of size $M+n$, which we may enumerate as $\ssb_1,\dots,\ssb_{M+n}$. Then for each $j$, incompressibility implies $\ssb_j \not\subseteq \bigsetcup_{k\colon k\neq j} \ssb_k$, so
$\bigcap_{k\colon k\neq j} \ssb_j\setminus \ssb_k$
is non-empty and we may pick some $\gm_j$ in this set.
Note that $\gm_1,\dots, \gm_{M+n}$ are all distinct.  Since at most $M$ of these belong to $\ssa$, at least $n$ of them belong to $\Om\setminus\ssa$; reordering if necessary, we may take these to be $\gm_1,\dots,\gm_n$.
\end{proof}

\begin{prop}\label{p:E and F}
  There is a sequence $(E_n)_{n=1}^\infty$ in $\FIN(\Om)$, and a sequence $(\cF_n)_{n=1}^\infty$ in $\cP(\cS)$, with the following properties:
 \begin{itemize}
 \item for each $n$, $|E_n|=|\cF_n|=n$;
   \item the sets $E_1,E_2,\dots$ are pairwise disjoint;
\item each $\ssx\in\cF_n$ satisfies $\ssx\supseteq E_j$ for all $1\leq j\leq n-1$;
  \item for each $n$ and each $\gm\in E_n$, there is a unique $\ssx\in\cF_n$ such that $\gm\in\ssx$.
  \end{itemize}
\end{prop}

\begin{proof}
We construct both sequences together by (strong) induction on $n$.
For the base case of $n=1$: pick any non-empty $\ssa\in\cS$, pick any $\gm\in \ssa$, and let$E_1=\{\gm\}$, $\cF_1=\{ \ssa\}$.

Now let $n\geq 2$ and suppose we have found $E_1,\dots, E_{n-1}$ and $\cF_1,\dots,\cF_{n-1}$ with the desired properties. In particular, since each point of $E_i$ belongs to some $\ssx\in\cF_i$, we have $E_i \subseteq \join(\cF_i)$ for $i=1,\dots, n-1$.
  
Let  $\ssa\defeq \join(\cF_1)\setcup \dots \setcup \join(\cF_{n-1})$; we have $\ssa\in\cS$ since $\cS$ is closed under finite unions. Now let $\ssb_1,\dots,\ssb_n\in\cS$ and $\gm_1,\dots,\gm_n\in\Om\setminus\ssa$ be as provided by Lemma~\ref{l:inductive step}.
For each $j$, let $\ssx_j\defeq \ssa\setcup \ssb_j\in\cS$ and note that $\ssx_j\supseteq\ssa\supseteq \bigsetcup_{i=1}^{n-1} E_i$. Also, since $\gm_j\notin \ssa$, whenever $j\neq k$ we have
\[ \gm_j \in ( \ssb_j\setminus \ssb_k )\setminus \ssa = \ssb_j \setminus \ssx_k \subseteq \ssx_j\setminus \ssx_k \;. \]
This shows that $\cF_n \defeq \{ \ssx_1,\dots, \ssx_n\}$ is incompressible (since $\gm_j$ is a point in $\ssx_j$ but not in $\bigsetcup_{k\colon k\neq j} \ssx_k$). Taking $E_n=\{\gm_1,\dots,\gm_n\}$, the inductive step is complete.
\end{proof}

\begin{dfn}[Defining our bad log-weight]
\label{d:bad log-weight}
  Keeping the notation of Proposition \ref{p:E and F}, let $D_n\defeq \bigsetcup_{i=1}^n E_i$ for each $n\in\Nat$; we also define $D_0\defeq\emptyset$ and $D_\infty=\bigsetcup_{n\geq 1} E_n$.
  Then, for each $\ssx\in\FIN(\Om)$, we make the following definitions:
  \begin{itemize}
  \item    $N(\ssx)\defeq \max \set{ n\in\Nat_0 \colon D_n\subseteq \ssx }$;
  \item $\red(\ssx) \defeq (\ssx\setcap D_\infty) \setminus D_{N(\ssx)} = (\ssx\setminus D_{N(\ssx)})\setcap D_\infty$;
    \item $\eta(x)\defeq |\red(\ssx)|$.
    \end{itemize}
\end{dfn}

\begin{prop}
The function $\eta:\FIN(\Om)\to [0,\infty)$ is subadditive.
\end{prop}

\begin{proof}
  Let $\ssx,\ssy\in\FIN(\Om)$. Since $N(\ssx\setcup\ssy)\geq\max(N(\ssx),N(\ssy))$, we have
  \[ D_{N(\ssx\setcup\ssy)} \supseteq D_{N(\ssx)}\setcup D_{N(\ssy)}. \]
 Hence
 $(    \ssx\setcup\ssy)\setminus D_{N(\ssx\setcup\ssy)}
    \subseteq ( \ssx\setminus D_{N(\ssx)}) \setcup  ( \ssy\setminus D_{N(\ssy)})$.
 Therefore, $\red(\ssx\cup\ssy)\subseteq \red(\ssx)\setcup\red(\ssy)$, and the result follows.
\end{proof}

\begin{lem}\label{l:barrier again}
Let $n\geq 2$ and let $\ssz_n\defeq \join(\cF_n)\in\cS$. Let $C\geq 0$ be such that $\ssz_n\in\FBP_C^\infty(\cF_n)$. Then $2C \geq |E_n|=n$.
\end{lem}

\begin{proof}
  Let $m\in\Nat$ be minimal with the following property: there exists some $\ssy\in \FBP_C^m(\cF_n)$ such that $\ssy\supseteq E_n$. This is well-defined, since $\ssz_n\supseteq E_n$ (as observed in the proof of Proposition \ref{p:E and F}) and $\ssz_n\in\FBP_C^\infty(\cF_n)$ by assumption.

  By definition of $\FBP_C$, there exist $\ssx_1,\ssx_2\in \FBP_C^{m-1}(\cF_n)$ such that $ \ssx_1\setcup \ssx_2\supseteq \ssy\supseteq E_n$.
  We claim that $\ssx_i\not\supseteq E_n$ for $i\in\{1,2\}$:
  \begin{itemize}
  \item[--] if $m\geq 2$, the claim follows from minimality of $m$;
  \item[--] if $m=1$, the claim holds because $\FBP_C^0(\cF_n)\subseteq \cF_n$ and for each $\ssx\in\cF_n$, $|\ssx\setcap E_n|=1$ while $|E_n|=n \geq 2$.
  \end{itemize}
  Therefore, $\red(\ssx_i)\supseteq E_n\setcap \ssx_i \neq \emptyset $, for $i\in\{1,2\}$, and so
  \[
  2C\geq \eta(\ssx_1)+\eta(\ssx_2) \geq |E_n\setcap \ssx_i| + |E_n\setcap \ssx_2| \geq | E_n\setcap (\ssx_1\setcup \ssx_2)| = |E_n|
  \]
as required.
\end{proof}

\begin{cor}
Let $\eta$ be as in Definition~\ref{d:bad log-weight}. Then $(S,\eta)$ does not have $1$-propagation.
\end{cor}

\begin{proof}
For $n\geq 2$, each $\ssx\in\cF_n$ satisfies $\ssx\setcap D_\infty = D_{n-1}\cup\{\gm\}$ for some $\gm\in E_n$. Therefore $\cF_n\subseteq  W_1$.
Now let $\ssz_n=\join(\cF_n)\in \filgen{\cF_n}$.
By the properties listed in
Proposition~\ref{p:E and F}, $\ssz_n\setcap D_\infty=D_n$, and $\eta(\ssz_n)=0$. Thus $\ssz_n\in  W_1$. On the other hand, by Lemma~\ref{l:barrier again},
$V_{\cF_n}(\ssz_n) \geq n/2$. Hence
\[
\sup_{\emptyset\neq \cE\subseteq W_1} \ \sup_{\ssz\in\filgen{\cE}\sscap  W_1} V_{\cE}(\ssz)
\geq \sup_{n\geq 2} V_{\cF_n}(\ssz_n) = + \infty
\]
and so $(S,\eta)$ does not propagate at level~$1$.
\end{proof}

\end{section}

\section*{Acknowledgements}
This paper forms part of a larger project, which grew out of conversations between the authors at the conference ``Banach Algebras and Applications'', held in Gothenburg, Sweden, July--August 2013, and which was further developed while the authors were attending the thematic program ``Abstract 
Harmonic Analysis, Banach and Operator Algebras''  at the Fields Institute, Canada, during March--April 2014. The authors thank the organizers of these meetings for invitations to attend and for pleasant environments to discuss research.

The first author acknowledges financial support to attend the latter meeting, in the form of a travel grant from the Faculty of Science and Technology at Lancaster University. The third author acknowledges the financial supports of a Fast Start Marsden Grant and of Victoria University of Wellington to attend both meetings.
Further work was done during the second author's visit to Lancaster University in November 2014, which was supported by a Scheme~2 grant from the London Mathematical Society
(reference 21313).
She thanks the Department of Mathematics and Statistics at Lancaster for their hospitality.

The first author thanks B. Horv\'ath and N. J. Lauststen for discussions in recent years concerning AMNM phenomena, and for their encouragement to write up the results presented here.
The second author also acknowledges support from National Science Foundation grant DMS-1902301 during the preparation of this article.


\newcommand{\etalchar}[1]{$^{#1}$}

\vfill

\newcommand{\address}[1]{{\small\sc#1.}}
\newcommand{\email}[1]{\texttt{#1}}

\noindent
\address{Yemon Choi,
Department of Mathematics and Statistics,
Lancaster University,
Lancaster LA1 4YF, United Kingdom} 

\email{y.choi1@lancaster.ac.uk}

\noindent
\address{Mahya Ghandehari,
Department of Mathematical Sciences,
University of Delaware,
Newark, Delaware 19716, United States of America}

\email{mahya@udel.edu}

\noindent
\address{Hung Le Pham,
School of Mathematics and Stat\-istics,
Victoria University of Wellington,
Wellington 6140, New Zealand}

\email{hung.pham@vuw.ac.nz}

\end{document}